\documentclass[a4paper,11pt]{amsart}
\usepackage{wrapfig}
\usepackage[dvips]{graphicx}
\usepackage{amsmath,amsthm,amssymb,amscd}
\theoremstyle{definition}
\newtheorem{thm}{Theorem}[section]

\newtheorem{lem}[thm]{Lemma}
\newtheorem{ex}[thm]{Example}
\newtheorem{rem}[thm]{Remark}
\theoremstyle{definition}

\begin{document}

\title[Trace scaling actions and fundamental groups]{A note on trace scaling actions and fundamental groups of $C^*$-algebras}
\author{Norio Nawata}
\address[Norio Nawata]{Graduate School of Mathematics, 
Kyushu University, Motooka, 
Fukuoka, 819-0395,  Japan}      
\email{n-nawata@math.kyushu-u.ac.jp}
\maketitle
\begin{abstract}
Using Effros-Handelman-Shen theorem and Elliott's classification theorem of AF algebras, 
we show that there exists a unital simple AF algebra $A$ with unique trace such that 
$A\otimes \mathbb{K}$ admits no trace scaling action of the fundamental group of $A$. 
\end{abstract}

\section{Introduction}
Let $M$ be a factor of type $\mathrm{II}_1$ with a normalized trace 
$\tau$. Murray and von Neumann introduced 
the fundamental group ${\mathcal F}(M)$ of $M$ in \cite{MN}. 
They showed that if $M$ is  hyperfinite, then 
${\mathcal F}(M) = {\mathbb R_+^{\times}}$. 
Since then 
there has been many works on the computation of 
the fundamental groups. Voiculescu \cite{Vo} showed that 
${\mathcal F}(L(\mathbb{F}_{\infty}))$ of the group factor 
of the free group $\mathbb{F}_{\infty}$ contains the positive rationals and 
Radulescu proved that 
${\mathcal F}(L(\mathbb{F}_{\infty})) = {\mathbb R}_+^{\times}$ in 
\cite{Ra}.  Connes \cite{Co} showed that if $G$ is an ICC group with property 
(T), then  ${\mathcal F}(L(G))$ is a countable group. Popa 
showed that any countable subgroup of $\mathbb R_+^{\times}$ 
can be realized as the fundamental group of some 
factor of type $\mathrm{II}_1$ in \cite{Po1}. 
Furthermore Popa and Vaes \cite{PV} exhibited a large family $\mathcal{S}$ 
of subgroups of $\mathbb{R}_{+}^\times$, containing $\mathbb{R}_{+}^\times$ 
itself, all of its countable subgroups, as well as uncountable subgroups with 
any Hausdorff dimension in $(0,1)$, such that for each $G\in\mathcal{S}$ 
there exist many free ergodic measure preserving actions of $\mathbb{F}_{\infty}$ 
for which the associated $\mathrm{II}_1$ factor $M$ has the fundamental group equal to $G$. 
In our previous paper \cite{NW} (see also \cite{N}), 
we introduced the fundamental group $\mathcal{F}(A)$ 
of a simple unital $C^*$-algebra $A$ with a normalized trace $\tau$ 
based on the computation  of Picard groups by 
Kodaka \cite{kod1}, \cite{kod2} and \cite{kod3}. 
The fundamental group ${\mathcal F}(A)$ is 
defined as the set of 
the numbers $\tau \otimes Tr(p)$ for some projection 
$p \in M_n(A)$ such that $pM_n(A)p$ is isomorphic to $A$. 
We computed the fundamental groups of several $C^*$-algebras and showed 
that any countable subgroup of $\mathbb{R}_+^\times$ 
can be realized as the fundamental group of a separable simple unital $C^*$-algebra 
with unique trace \cite{NW2}. 

The fundamental group of a $\mathrm{II}_1$ factor $M$ is equal to the set of trace-scaling constants 
for automorphisms of $M\otimes B(\mathcal{H})$. 
We have a similar fact, that is, the fundamental group 
of a $C^*$-algebra $A$ is equal to the set of trace-scaling constants for 
automorphisms of $A\otimes \mathbb{K}$ \cite{NW} (see also \cite{N}). 
It is of interest to know whether $A\otimes\mathbb{K}$ admits a trace scaling action of 
$\mathcal{F}(A)$. 
In the case where $M$ is a factor of type $\mathrm{II}_1$, the existence of a trace scaling 
(continuous) $\mathbb{R}_{+}^\times$-action on $M\otimes B(\mathcal{H})$ is equivalent to 
the existence of a type $\mathrm{III}_1$ factor having a core isomorphic to $M\otimes B(\mathcal{H})$ 
by the continuous decomposition of type $\mathrm{III}_1$ factors. (See \cite{T} and \cite{CT}.) 
Hence this question is important in the theory of von Neumann algebras. 
Radulescu showed that $L(\mathbb{F}_{\infty})\otimes B(\mathcal{H})$ admits a trace scaling 
action of $\mathbb{R}_{+}^\times$ in \cite{Ra2}. Therefore there exists a type $\mathrm{III}_1$ 
factor having a core isomorphic to $L(\mathbb{F}_{\infty})\otimes B(\mathcal{H})$. 
Popa and Vaes \cite{PV2} showed that there exists a $\mathrm{II}_1$ factor $M$ such that 
$\mathcal{F}(M)=\mathbb{R}_{+}^\times$ and $M\otimes B(\mathcal{H})$ admits no 
trace scaling (continuous) action of $\mathbb{R}_{+}^\times$. 

In this paper we consider trace scaling actions on certain AF algebras. 
If $A$ is a UHF algebra, then $A\otimes\mathbb{K}$ admits a trace scaling 
action of $\mathcal{F}(A)$. 
Using Effros-Handelman-Shen theorem and Elliott's classification theorem 
of AF algebras, we show that there exists a unital simple AF algebra $A$ with unique trace 
such that $A\otimes \mathbb{K}$ admits no trace scaling action of $\mathcal{F}(A)$. 
Note that there exist remarkable works of the classification of trace scaling automorphisms in 
\cite{BK}, \cite{EEK} and \cite{EK}. 
But we do not consider the classification of trace scaling actions in this paper.

\section{Examples}\label{sec:exp} 
We recall some definitions in \cite{NW}. Let $A$ be a unital simple $C^*$-algebra with a unique 
normalized trace $\tau$ and $Tr$ the usual unnormalized trace on $M_n(\mathbb{C})$. 
Put 
$$
\mathcal{F}(A) :=\{ \tau\otimes Tr(p) \in \mathbb{R}^{\times}_{+}\ | \ 
 p \text{ is a projection in } M_n(A) \text{ such that } pM_n(A)p  \cong A \}. 
$$
Then $\mathcal{F}(A)$ is a multiplicative subgroup of $\mathbb{R}_{+}^\times$ 
by Theorem 3.1 in \cite{NW}. 
For an additive subgroup $E$ of 
$\mathbb{R}$ containing 1, 
we define the positive inner multiplier group $IM_+(E)$ of $E$  by 
$$
IM_+(E) = \{t \in {\mathbb R}_+^{\times} \ | t \in E, t^{-1}  \in E, \text{ and } 
        tE = E \}. 
$$
Then we have $\mathcal{F}(A) \subset IM_+(\tau_*(K_0(A)))$ by Proposition 3.7 in \cite{NW}. 
This obstruction enables us to compute fundamental groups easily. 
For $x\in (A\otimes \mathbb{K})_+$, set 
$\hat{\tau} (x)=\sup\{\tau \otimes Tr(y):y\in \cup _n M_n(A), y\leq x\}$. 
Define $\mathcal{M}_\tau^+=
\{x\geq 0: \hat{\tau}  (x)<\infty \}$ and $\mathcal{M}_\tau =\mathrm{span} \mathcal{M}_\tau^+$. 
Then $\hat{\tau}$ is a densely defined (with the domain $\mathcal{M}_\tau $) 
lower semicontinuous trace 
on $A\otimes \mathbb{K}$. 
Since the normalize trace on a unital $C^*$-algebra $A$ is unique, the lower semicontinuous 
densely defined trace on $A\otimes \mathbb{K}$ is unique up to constant multiple. 
It is clear that for any $\alpha\in \mathrm{Aut}(A\otimes \mathbb{K})$, 
$\hat{\tau}\circ \alpha$ is a densely defined 
(with the domain $\alpha^{-1}(\mathcal{M}_\tau )$) lower semicontinuous trace on 
$A\otimes\mathbb{K}$. 
Therefore there exists a positive number $\lambda$ such that 
$\hat{\tau}\circ \alpha =\lambda\hat{\tau}$, and hence 
$\alpha^{-1}(\mathcal{M}_\tau )=\mathcal{M}_\tau$. 
We define the set of trace-scaling constants for automorphisms:  
$$
\mathfrak{S}(A)
:= \{ \lambda \in \mathbb{R}^{\times}_+ \ | \ 
\hat{\tau} \circ \alpha = \lambda \hat{\tau} \text{ for some  } 
 \alpha \in \mathrm{Aut} (A\otimes \mathbb{K}) \ \}. 
$$
Then $\mathcal{F}(A)=\mathfrak{S}(A)$ by Proposition 3.28 in \cite{NW}. 
Therefore it is of interest to know whether $A\otimes\mathbb{K}$ admits a trace scaling action of 
$\mathcal{F}(A)$. 

It is clear that if the fundamental group of $A$ is singly generated, $A\otimes\mathbb{K}$ 
admits a trace scaling action of $\mathcal{F}(A)$. See \cite{NW} and \cite{NW2} for such 
examples. We shall show some examples of AF algebras $A$ such that $A\otimes\mathbb{K}$ 
admits a trace scaling action of $\mathcal{F}(A)$. 
\begin{ex}\label{ex:uhf}
Consider a UHF algebra $M_{2^\infty 3^\infty}$. 
Then the fundamental group of $M_{2^\infty 3^\infty}$ is a multiplicative subgroup generated by 
$2$ and $3$. Hence $\mathcal{F}(M_{2^\infty 3^\infty})$ is isomorphic to $\mathbb{Z}^2$ as a 
group. Since $M_{2^\infty 3^\infty}\otimes\mathbb{K}$ is isomorphic to 
$M_{2^\infty}\otimes\mathbb{K}\otimes M_{3^\infty}\otimes\mathbb{K}$, there exists 
a trace scaling $\mathbb{Z}^2$-action on $M_{2^\infty 3^\infty}\otimes\mathbb{K}$. 
In general, if $A$ is a UHF algebra, then $\mathcal{F}(A)$ is a free abelian group 
(see \cite{NW}) and $A\otimes\mathbb{K}$ admits a trace scaling action of $\mathcal{F}(A)$. 
\end{ex}
\begin{ex}\label{ex:admit} 
Let $A$ be a unital simple AF algebra such that 
$K_0(A)=\mathbb{Z}+\mathbb{Z}\sqrt{3}$, $K_0(A)_{+}=(\mathbb{Z}+\mathbb{Z}\sqrt{3})\cap\mathbb{R}_{+}$ 
and $[1]_0=1$. Then $\mathcal{F}(A)=\{(2+\sqrt{3})^n:n\in\mathbb{Z}\}$ 
(see Proposition 3.17 and Corollary 3.18 in \cite{NW}). 
Consider $B=M_{5^\infty}\otimes A$. Then it is easily seen that $\tau_{*}(K_0(B))= 
\mathbb{Z}[\frac{1}{5}]+\mathbb{Z}[\frac{1}{5}]\sqrt{3}$ and 
$\tau_{*}$ is an order isomorphism. 
We shall show that $IM_+(\tau_*(K_0(B)))$ is generated by $5$ and $2+\sqrt{3}$. 
Since $\tau_{*}(K_0(B))$ is a subring of $\mathbb{R}$, $IM_+(\tau_*(K_0(B)))$ is a group of 
positive invertible elements. 
Define a multiplicative map $N$ of $\mathbb{Z}[\frac{1}{5}]+\mathbb{Z}[\frac{1}{5}]\sqrt{3}$ 
to $\mathbb{Z}[\frac{1}{5}]$ by $N(a+b\sqrt{3})=a^2-3b^2$ for any 
$a,b\in\mathbb{Z}[\frac{1}{5}]$. If $a+b\sqrt{3}$ is an invertible element in 
$\mathbb{Z}[\frac{1}{5}]+\mathbb{Z}[\frac{1}{5}]\sqrt{3}$, then there exists an integer $n$ 
such that $N(a+b\sqrt{3})=\pm 5^n$. 
Elementary computations shows that 
$x^2-3y^2 \equiv 0\; \mathrm{mod}\; 25$ implies $x\equiv 0\;\mathrm{mod}\;5$ and 
$y\equiv 0\;\mathrm{mod}\; 5$. 
It is easy to see that no integers $x$ and $y$ satisfy equations $x^2-3y^2 =\pm 5$ 
or $x^2-3y^2=-1$. There exist integers $x$ and $y$ satisfy the equation $x^2-3y^2=-1$. 
(See, for example, \cite{JW} and \cite{NW}.) 
Therefore it can be easily checked that $IM_+(\tau_*(K_0(B)))$ is generated by $5$ and 
$2+\sqrt{3}$. 
Hence we see that 
$\mathcal{F}(B)=\{5^n(2+\sqrt{3})^m:n,m\in\mathbb{Z} \}$ by Proposition 3.18 in \cite{NW} 
and $B\otimes \mathbb{K}$ admits a trace scaling action. 
\end{ex}
We shall show that there exists a unital simple AF algebra $A$ with unique trace 
such that $A\otimes \mathbb{K}$ admits no trace scaling action of $\mathcal{F}(A)$. 
Define 
$$
E=\{(\frac{j+k\sqrt{3}}{5^{6i}},
\left (\begin{array}{c} 
           x \\
           y  \end{array} \right )) \in \mathbb{R} \times \mathbb{Z}^{2} 
\ | \ i,j,k,x,y\in\mathbb{Z},x\equiv j\; \mathrm{mod}\; 9,y\equiv k\; \mathrm{mod}\; 3 \}
$$ 
$$
E_+=\{(r,\left (\begin{array}{c} 
           x \\
           y  \end{array} \right ))\in E:r>0\}\cup \{(0,\left(\begin{array}{c} 
           0 \\
           0  \end{array} \right ))\}
\ \ \text{and} \ \  [u]_0= (1,\left(\begin{array}{c} 
           1 \\
           0  \end{array} \right )). 
$$
Then there exists a simple $AF$-algebra $A$ with a unique normalized trace $\tau$ 
such that $(K_0(A),K_0(A)_+,[1_A]_0)=(E,E_+,u)$ by Effros-Handelman-Shen theorem \cite{EHS}. 
\begin{lem}\label{lem:group}
With notation as above the fundamental group of $A$ is equal to the multiplicative group 
generated by $5$ and $2+\sqrt{3}$. 
\end{lem}
\begin{proof}
Since $\tau_*(K_0(A))$ is equal to $\mathbb{Z}[\frac{1}{5}]+\mathbb{Z}[\frac{1}{5}]\sqrt{3}$, 
$\mathcal{F}(A)$ is a subgroup of $\{5^n(2+\sqrt{3})^m:n,m\in\mathbb{Z} \}$ by an argument 
in Example \ref{ex:admit}. 
Define an additive homomorphism $\phi :E\rightarrow E$ by 
$$\phi ((r,\left (\begin{array}{c} 
           x \\
           y  \end{array} \right ))) =
           (5r,\left (\begin{array}{cc}
           5 & 9 \\
           6 & 11 \end{array} \right ) \left (\begin{array}{c} 
           x \\
           y  \end{array} \right )).
$$          
Computations show that $\phi$ is a well-defined order isomorphism of $E$ with 
$\phi (u)=(5,\left (\begin{array}{c} 
           5 \\
           6  \end{array} \right ))$. There exist a natural number $n$ and a projection $p$ 
in $M_n(A)$ such that $[p]_0=(5, \left (\begin{array}{c}  
           5 \\
           6  \end{array} \right ))$ and $\tau\otimes Tr (p)=5$. 
Since we have 
$(K_0(pM_n(A)p),K_0(pM_n(A)p)_+,[p]_0) =
(E,E_+,(5, \left (\begin{array}{c} 
           5 \\
           6  \end{array} \right )))$, there exists an isomorphism $f:A\rightarrow pM_n(A)p$ 
with $f_* = \phi$ by Elliott's classification theorem of AF algebra \cite{E}. 
Therefore $5\in\mathcal{F}(A)$. Define an additive homomorphism $\psi :E\rightarrow E$ by 
$$\psi ((r,\left (\begin{array}{c} 
           x \\
           y  \end{array} \right ))) =
           ((2+\sqrt{3})r,\left (\begin{array}{cc}
           2 & 3 \\
           1 & 2 \end{array} \right ) \left (\begin{array}{c} 
           x \\
           y  \end{array} \right )).
$$          
Then we see that $2+\sqrt{3}\in\mathcal{F}(A)$. 
Consequently $\mathcal{F}(A)$ is the multiplicative group generated by $5$ and $2+\sqrt{3}$. 
\end{proof}
We shall consider the order automorphisms of $(E,E_{+})$. 
\begin{lem}\label{lem:commute}
Let $\phi$ be an order automorphism of $(E,E_{+})$. Then there exist integers $a,b,c,d$ and 
a positive invertible element $\lambda$ in $\mathbb{Z}[\frac{1}{5}]+\mathbb{Z}[\frac{1}{5}]\sqrt{3}$ 
such that  $ad-bc=\pm 1$ and 
$$\phi ((r,\left (\begin{array}{c} 
           x \\
           y  \end{array} \right )))=
           (\lambda r,\left (\begin{array}{cc}
           a & b \\
           c & d \end{array} \right ) \left (\begin{array}{c} 
           x \\
           y  \end{array} \right ) ).$$ 
Moreover if $\lambda =5$, then 
$$
           \left (\begin{array}{cc}
           a & b \\
           c & d \end{array} \right )\equiv 
\left (\begin{array}{cc}
           5 & 0 \\
           0 & 2 \end{array} \right ), 
\left (\begin{array}{cc}
           5 & 0 \\
           3 & 2 \end{array} \right ),
\left (\begin{array}{cc}
           5 & 0 \\
           6 & 2 \end{array} \right )
\; \mathrm{mod}\; 9
$$ 
and if $\lambda =2+\sqrt{3}$, then 
$$
           \left (\begin{array}{cc}
           a & b \\
           c & d \end{array} \right )\equiv 
\left (\begin{array}{cc}
           2 & 3 \\
           1 & 2 \end{array} \right ), 
\left (\begin{array}{cc}
           2 & 3 \\
           4 & 2 \end{array} \right ),
\left (\begin{array}{cc}
           2 & 3 \\
           7 & 2 \end{array} \right )
\; \mathrm{mod}\; 9.
$$ 
\end{lem}
\begin{proof}
We denote by $(\phi_1((r,\left (\begin{array}{c} 
           x \\
           y  \end{array} \right ))), \phi_2((r,\left (\begin{array}{c} 
           x \\
           y  \end{array} \right ))))$ 
           the element $\phi((r,\left (\begin{array}{c} 
           x \\
           y  \end{array} \right )))$ for any 
           $(r,\left (\begin{array}{c} 
           x \\
           y  \end{array} \right ))\in E$. 
Consider a subgroup $F$ generated by $(0,\left (\begin{array}{c} 
           9 \\
           0  \end{array} \right ))$ and 
           $(0,\left (\begin{array}{c} 
           0 \\
           3  \end{array} \right ))$. 
Then $F$ is an $\phi$-invariant subgroup because $\phi$ is an order isomorphism. 
Hence there exist integers $m_1$, $m_2$, $m_3$ and $m_4$ such that 
$m_1m_4-m_2m_3=\pm 1$ and 
$\phi_2 ((0,\left (\begin{array}{c} 
           x \\
           y  \end{array} \right )))= 
           \left (\begin{array}{cc}
           m_1 & 3m_2 \\
           \frac{m_3}{3} & m_4 \end{array} \right )
           \left (\begin{array}{c} 
           x \\
           y  \end{array} \right )$ for any $(0,\left (\begin{array}{c} 
           x \\
           y  \end{array} \right ))\in F$. 
Furthermore we see that there exists a positive invertible element $\lambda$ in 
$\mathbb{Z}[\frac{1}{5}]+\mathbb{Z}[\frac{1}{5}]\sqrt{3}$ such that 
$\phi_1 ((r,\left (\begin{array}{c} 
           x \\
           y  \end{array} \right )))=\lambda r
           $. 
Since $5^{6i}\phi ((\frac{9}{5^{6i}},\left (\begin{array}{c} 
           0 \\
           0  \end{array} \right )))= \phi ((9,\left (\begin{array}{c} 
           0 \\
           0  \end{array} \right )))$ for any $i\in\mathbb{Z}$, 
we see that $\phi ((9,\left (\begin{array}{c} 
           0 \\
           0  \end{array} \right )))= (9\lambda ,\left (\begin{array}{c} 
           0 \\
           0  \end{array} \right ))$. 
This observation and easy computations show that 
$\phi ((1,\left (\begin{array}{c} 
           1 \\
           0  \end{array} \right )))=
(\lambda ,\left (\begin{array}{cc}
           m_1 & 3m_2 \\
           \frac{m_3}{3} & m_4 \end{array} \right )
           \left (\begin{array}{c} 
           1 \\
           0  \end{array} \right ))$ and $\frac{m_3}{3}\in\mathbb{Z}$. 
In a similar way, we see that 
$\phi ((\sqrt{3},\left (\begin{array}{c} 
           0 \\
           1  \end{array} \right )))=
(\lambda ,\left (\begin{array}{cc}
           m_1 & 3m_2 \\
           \frac{m_3}{3} & m_4 \end{array} \right )
           \left (\begin{array}{c} 
           0 \\
           1  \end{array} \right ))$. 
It is easily seen that $\phi$ is determined by the values of 
$\phi ((1,\left (\begin{array}{c} 
           1 \\
           0  \end{array} \right )))$, 
$\phi ((\sqrt{3},\left (\begin{array}{c} 
           0 \\
           1  \end{array} \right )))$, 
$\phi ((0,\left (\begin{array}{c} 
           9 \\
           0  \end{array} \right )))$ and 
$\phi ((0,\left (\begin{array}{c} 
           0 \\
           3  \end{array} \right )))$. 
Therefore there exist integers $a,b,c,d$ and 
a positive invertible element $\lambda$ in $\mathbb{Z}[\frac{1}{5}]+\mathbb{Z}[\frac{1}{5}]\sqrt{3}$ 
such that  $ad-bc=\pm 1$ and 
$$\phi ((r,\left (\begin{array}{c} 
           x \\
           y  \end{array} \right )))=
           (\lambda r,\left (\begin{array}{cc}
           a & b \\
           c & d \end{array} \right ) \left (\begin{array}{c} 
           x \\
           y  \end{array} \right ) ).$$ 

Let $\lambda =5$, then $a\equiv 5\;\mathrm{mod}\; 9$, 
$b\equiv 0\;\mathrm{mod}\; 9$, $c\equiv 0\;\mathrm{mod}\; 3$ and 
$d\equiv 5\;\mathrm{mod}\; 3$ by the definition of $E$. 
If $ad-bc=1$, then $d\equiv 5^5\;\mathrm{mod}\; 9$, $-b\equiv 0\;\mathrm{mod}\; 9$, 
$-c\equiv 0\;\mathrm{mod}\; 3$ and $a\equiv 5^5\;\mathrm{mod}\; 3$ because $\phi$ is 
an isomorphism. Therefore computations show 
$$
           \left (\begin{array}{cc}
           a & b \\
           c & d \end{array} \right )\equiv 
\left (\begin{array}{cc}
           5 & 0 \\
           0 & 2 \end{array} \right ), 
\left (\begin{array}{cc}
           5 & 0 \\
           3 & 2 \end{array} \right ),
\left (\begin{array}{cc}
           5 & 0 \\
           6 & 2 \end{array} \right )
\; \mathrm{mod}\; 9
.$$ 
If $ad-bc=-1$, then then $-d\equiv 5^5\;\mathrm{mod}\; 9$, $b\equiv 0\;\mathrm{mod}\; 9$, 
$c\equiv 0\;\mathrm{mod}\; 3$ and $-a\equiv 5^5\;\mathrm{mod}\; 3$. 
There does not exist a integer $a$ such that $a\equiv 5\;\mathrm{mod}\; 9$ 
and $-a\equiv 5^5\;\mathrm{mod}\; 3$. Therefore we reach a conclusion in the case 
$\lambda =5$. 
In the case $\lambda =2+\sqrt{3}$, the similar argument as above proves the lemma. 
\end{proof}
\begin{thm}\label{thm:scaling main}
There exists a unital simple AF algebra $A$ with unique trace 
such that $A\otimes \mathbb{K}$ admits no trace scaling action of $\mathcal{F}(A)$. 
\end{thm}
\begin{proof}
Let 
$$
E=\{(\frac{j+k\sqrt{3}}{5^{6i}},
\left (\begin{array}{c} 
           x \\
           y  \end{array} \right )) \in \mathbb{R} \times \mathbb{Z}^{2} 
\ | \ i,j,k,x,y\in\mathbb{Z},x\equiv j\; \mathrm{mod}\; 9,y\equiv k\; \mathrm{mod}\; 3 \}
$$ 
$$
E_+=\{(r,\left (\begin{array}{c} 
           x \\
           y  \end{array} \right ))\in E:r>0\}\cup \{(0,\left(\begin{array}{c} 
           0 \\
           0  \end{array} \right ))\}
\ \ \text{and} \ \  [u]_0= (1,\left(\begin{array}{c} 
           1 \\
           0  \end{array} \right )). 
$$
Then there exists a simple AF algebra $A$ with a unique normalized trace $\tau$ 
such that $(K_0(A),K_0(A)_+,[1_A]_0)=(E,E_+,u)$ by Effros-Handelman-Shen theorem \cite{EHS}. 
By Lemma \ref{lem:group}, $\mathcal{F}(A)=\{5^n(2+\sqrt{3})^m:n,m\in\mathbb{Z} \}$. 
Let $\alpha$ be an automorphism of $A\otimes\mathbb{K}$ such that 
$\hat{\tau}\circ \alpha =5\hat{\tau}$ and $\beta$ an automorphism of $A\otimes\mathbb{K}$ 
such that $\hat{\tau}\circ \beta =(2+\sqrt{3})\hat{\tau}$. 
Then $\alpha_{*}$ and $\beta_{*}$ are order isomorphisms of $(K_0(A),K_0(A)_+)$. 
Lemma \ref{lem:commute} and computations show that 
$\alpha_{*}\circ\beta_{*}\neq \beta_{*}\circ \alpha_{*}$. 
Therefore $A\otimes \mathbb{K}$ admits no trace scaling action of $\mathcal{F}(A)$. 
\end{proof}
\begin{rem}
Let $A$ be a unital simple $C^*$-algebra with a unique normalized trace $\tau$. 
We denote by $\mathrm{Pic}(A)$ the Picard group of $A$ (see \cite{BGR}). 
Assume that the normalized trace on $A$ separates equivalence classes of projections. 
Then we have the 
following exact sequence \cite{NW} (see also \cite{kod1}). 
\[\begin{CD}
      {1} @>>> \mathrm{Out}(A) @>\rho_A>> \mathrm{Pic}(A) @>T>> \mathcal{F}(A)
 @>>> {1} \end{CD}. \] 
If $A\otimes\mathbb{K}$ admits a trace scaling action of $\mathcal{F}(A)$, then 
$\mathrm{Pic}(A)$ is isomorphic to a semidirect product of $\mathrm{Out}(A)$ with 
$\mathcal{F}(A)$. Example \ref{ex:uhf} and Example \ref{ex:admit} are such examples. 
We do not know whether there exists a simple $C^*$-algebra $A$ with a unique normalized 
trace $\tau$ such that the normalized trace on $A$ separates 
equivalence classes of projections and $A\otimes\mathbb{K}$ admits no trace scaling 
action of $\mathcal{F}(A)$. 
\end{rem}
\begin{rem}
If $A$ is a $C^*$-algebra in the proof of Theorem \ref{thm:scaling main}, then 
it can be checked that $\mathrm{Out}(A)$ is not a normal subgroup of $\mathrm{Pic}(A)$ 
by Lemma \ref{lem:commute}, Proposition 1.5 in \cite{kod1} and Elliott's classification 
theorem of AF algebras. 
\end{rem}

\end{document}